\documentclass[a4paper,12pt]{amsart}

\usepackage{amssymb}
\usepackage{amsthm}
\usepackage{color}
\input xy
\xyoption{all}
\usepackage{hyperref}
\hypersetup{colorlinks, linkcolor=blue, filecolor=green, urlcolor=blue, citecolor=blue}

\theoremstyle{plain}
\newtheorem{theorem}{Theorem}[section]
\newtheorem{lemma}[theorem]{Lemma}
\newtheorem{corollary}[theorem]{Corollary}
\newtheorem{proposition}[theorem]{Proposition}
\theoremstyle{definition}
\newtheorem{definition}[theorem]{Definition}
\newtheorem{example}[theorem]{Example}
\newtheorem{remark}[theorem]{Remark}
\numberwithin{equation}{section}

\newcommand{\Z}{\mathbb{Z}}
\newcommand{\N}{\mathbb{N}}

\newcommand{\im}{\operatorname{Im}}
\newcommand{\Path}{\operatorname{Path}}

\newcommand{\Ext}{(\widehat{E})^0}

\begin{document}
\title[Quotients of LPAs over rings by graded ideals]{Quotients of Leavitt path algebras over commutative rings by graded ideals}
\author{Sevgi Harman, Müge Kanuni, Guillermo Vera de Salas}
\address{Department of Mathematics Engineering, İstanbul Technical University, 34467 İstanbul, Turkey.}
\email{harman@itu.edu.tr}
\address{Department of Mathematics Engineering, Özyeğin University,
İstanbul, 34794, Turkey.}
\email{muge.er@ozyegin.edu.tr}
\address{ Departamento de Matem\'{a}tica  Aplicada, Ciencia e Ingenier\'{\i}a de los Materiales y Tecnolog\'{\i}a Electr\'onica,
Universidad Rey Juan Carlos, 28933 M\'{o}s\-to\-les (Madrid), Spain}
\email{guillermo.vera@urjc.es}

\keywords{Leavitt path algebra, graded ideal, $I$-basic graded ideal, graded ideal function, quotient of Leavitt path algebra}
\subjclass[2020]{16S88, 16D25, 16W50}
\date{\today}

\begin{abstract}
Let $E$ be a graph and $R$ a unital commutative ring. We study quotients of the Leavitt path algebra $L_R(E)$ by graded ideals which, unlike in the case of field coefficients, need not be basic. Using the saturated functions of Rigby and van den Hove, we introduce $I$-basic graded ideals, where $I$ is an ideal of $R$, and we prove that, for an arbitrary graph $E$, the quotient of $L_R(E)$ by an $I$-basic graded ideal is isomorphic to the Leavitt path algebra of a quotient graph of $E$ with coefficients in $R/I$. More precisely, for every graded ideal there is a canonical epimorphism from such an algebra onto the quotient, and it is an isomorphism if and only if the graded ideal is $I$-basic. This extends a result of Larki for basic graded ideals. We also introduce graded ideal functions, defined on the vertices of $E$ and on the pairs $(H,v)$ formed by a hereditary saturated set $H$ and a breaking vertex $v$ of $H$, and we show that they form a lattice isomorphic to the lattice of graded ideals of $L_R(E)$. As an application, for a finite graph $E$ in which no cycle has an exit, we describe the quotient of $L_R(E)$ by any graded ideal as a direct sum, indexed by the sinks and the cycles of $E$, of Leavitt path algebras over quotient rings of $R$, namely matrix algebras over quotient rings of $R$ and of $R[x,x^{-1}]$. Finally, examples on the infinite binary tree and the Toeplitz graph show that the hypothesis on the finiteness and cycles cannot be removed.
\end{abstract}

\maketitle

\section{Introduction}
Although Leavitt path algebras have been a fruitful avenue of research in the last two decades, they have mainly been studied over a field. Leavitt path algebras over unital commutative rings were studied by M. Tomforde \cite{Tomforde} and H. Larki \cite{Larki}. In order to understand the ideal structure of Leavitt path algebras over unital commutative rings, the notion of a basic ideal was introduced. When the ring is a field, all ideals are basic; however, there are non-basic ideals in the general setting. M. Tomforde \cite[Theorem 7.9]{Tomforde} gave a one-to-one correspondence between the hereditary saturated subsets of a row-finite graph and the basic graded ideals of the Leavitt path algebra. Later, H. Larki extended the definition of a basic ideal to arbitrary graphs and obtained a similar result for admissible pairs \cite[Theorem 3.10]{Larki}.
A paper on the ideals of Steinberg algebras by L.O. Clark, C. Edie-Michell, A. an Huef and A. Sims \cite{IdealsSteinberg} gives the characterization of both basic and non-basic graded ideals of a Leavitt path algebra over a row-finite graph with no sinks and satisfying Condition (K). S. Rigby and T. van den Hove \cite{Rigby} then studied basic, non-basic, graded and non-graded ideals of a Leavitt path algebra over a unital commutative ring, and gave a complete characterization of the two-sided ideals of a Leavitt path algebra of an arbitrary graph with coefficients in a unital commutative ring. For this characterization, \cite{Rigby} uses a function $f$, called a saturated function, from the set of admissible pairs of the graph into the lattice of ideals of the ring, and \cite{IdealsSteinberg} uses analogous maps defined on hereditary saturated sets. Our approach in this paper uses both the saturated function $f$ and a function $\varphi$, defined on a set $\Ext$ of extended vertices, with values in the lattice of ideals of the ring. The set $\Ext$ consists of the vertices of $E$ together with one symbol $v^H$ for each pair $(H,v)$, where $H$ is a hereditary saturated set and $v$ is a breaking vertex of $H$.

We now explain the role of the function $\varphi$. For a basic graded ideal $I(H,S)$, the quotient $L_R(E)/I(H,S)$ is the Leavitt path algebra of the quotient graph $E\backslash(H,S)$, obtained by removing the vertices of $H$ and adding a sink $v'$ for each $v \in B_H \backslash S$. This fails for non-basic graded ideals: a vertex may lie in a graded ideal $A$ only after being multiplied by the elements of some ideal of $R$, and every vertex may appear in $A$ with non-zero coefficients, so that no vertex is removed. We therefore consider, for each vertex $u$, the ideal $\varphi(u)$ of the coefficients $k$ with $ku \in A$ and, for each pair $(H,v)$ with $v \in B_H$, the ideal $\varphi(v^H)$ of the $k$ with $kv^H \in A$ and $kH \subseteq A$. These are the coefficients of the generators of $A$ (Theorem \ref{th:iso}(2); compare \cite[Theorem 4]{Ranga} over a field and \cite[Corollary 5.6]{Rigby} over a commutative ring). The values of $\varphi$ are not independent: $\varphi(s(e)) \subseteq \varphi(r(e))$ for every edge $e$, and at a regular vertex $u$ the ideal $\varphi(u)$ is the intersection of the values at the ranges of its edges; for row-finite graphs, these conditions characterize graded ideal functions (Proposition \ref{prop:local}). Thus the coefficients propagate backwards along the paths of $E$. When $E$ is finite and no cycle has an exit, every vertex reaches a sink or a cycle, and the quotient $L_R(E)/A$ is determined by the values of $\varphi$ at the sinks and at the vertices on cycles (Theorem \ref{Th:MainI}).

In Section \ref{Sec:gradedidealfunction}, we introduce graded ideal functions (Definition \ref{Def:varphi}) and prove that the lattice of graded ideal functions is isomorphic to the lattice of saturated functions and to the lattice of graded ideals of $L_R(E)$ (Theorem \ref{th:iso}). For row-finite graphs, graded ideal functions are characterized by a local condition at each vertex (Proposition \ref{prop:local}).

In Section \ref{Sec:Main_via_f}, the quotient of a Leavitt path algebra by an $I$-basic graded ideal (Definition \ref{Def:I-basicgradedideal}) is described by using the saturated function $f$. The $I$-basic graded ideals are precisely the ideals $I(H,S) + IL_R(E)$ (Lemma \ref{lem:Ibasic}). The main result is Theorem \ref{cor:main}: if $E$ is an arbitrary graph and $A$ is a graded ideal of $L_R(E)$, then the canonical epimorphism $L_{R/I}(E\backslash(H,S)) \rightarrow L_R(E)/A$ is an isomorphism if and only if $A$ is $I$-basic, where $H = A \cap E^0$, $S$ is the set of breaking vertices $v$ of $H$ with $v^H \in A$, and $I$ is the smallest value of the saturated function of $A$. This extends the result of H. Larki \cite[Theorem 3.10]{Larki}. Examples are given to illustrate the result on both row-finite graphs and graphs with infinite emitters.

In Section \ref{Sec:Main_via_varphi}, we use graded ideal functions to study quotients by arbitrary (basic or non-basic) graded ideals. For a finite graph $E$ in which no cycle has an exit, Theorem \ref{Th:MainI} shows that
$$L_R(E)/A \cong \bigoplus_{I \neq R} L_{R/I}\left(_{X_I}E\right) \cong \bigoplus_{w} M_{n(w)}\left(R/\varphi(w)\right) \oplus \bigoplus_{C} M_{m(C)}\big((R/\varphi(C^0))[x,x^{-1}]\big),$$
where $I$ runs over the values different from $R$ taken by $\varphi$ at the sinks and at the vertices on cycles, $X_I$ is the set of these vertices at which $\varphi$ takes the value $I$, $w$ runs over the sinks, $C$ runs over the cycles up to equivalence, and $n(w)$, $m(C)$ are suitable numbers of paths. In this setting, the sinks and the vertices on cycles are the terminal vertices of L. Va\v{s} \cite{VasPQ}, and the new vertices of the graph $_{X_I}E$ (Definition \ref{Gporcupine}) correspond to the paths that meet $X_I$ only at their last vertex. This relies on the decomposition of $L_R(E)$ into matrix algebras, obtained by G. Abrams, G. Aranda Pino and M. Siles Molina \cite{AAPSM, AAPSM2} over fields and by A. Ko\c{c} and M. Özaydın \cite[Proposition 15]{KO} over commutative rings. The graph $_XE$  originated from the {\it jellyfish graph}, presented by M. Özaydın in lectures at the CIMPA research school ``Leavitt path algebras and graph C*-algebras'' (Turkey, 2015), and was independently discovered and published by L. Va\v{s} as the {\it porcupine graph} \cite[Definition 3.1]{Vas}; in Definition \ref{Gporcupine}, $X$ need not be hereditary or saturated. Example \ref{ex:tree}, the infinite binary tree, shows that the finiteness of $E$ is needed, even for row-finite graphs, and we conclude with Example \ref{ex:Toeplitz2}, which shows that the hypothesis on the cycles cannot be removed: for the Toeplitz graph $T$ and the graded ideal $A$ of $L_{\Z}(T)$ generated by $2v$ and $4u$, the analogue of Theorem \ref{Th:MainI} would predict $L_{\Z}(T)/A \cong M_\infty(\Z_2) \oplus \Z_4[x,x^{-1}]$, but $L_{\Z}(T)/A$ is indecomposable.

\section{Preliminaries}
Throughout the paper, $R$ denotes a unital commutative ring, and $\Z_n = \Z/n\Z$.

\subsection{Leavitt path algebras}
A \textit{directed graph} $E = (E^0, E^1,r,s)$ consists of a set $E^0$ of vertices, a set $E^1$ of edges and the range and source maps $r, s \colon E^1 \rightarrow E^0$. A vertex $v$ is a \textit{sink} if $s^{-1}(v) = \emptyset$, \textit{regular} if $0 < |s^{-1}(v)| < \infty$, and an \textit{infinite emitter} if $|s^{-1}(v)|=\infty$. A \textit{path} $p = e_1 \cdots e_n$ ($n \geq 1$) is a sequence of edges with $r(e_i) = s(e_{i+1})$; we write $|p| = n$, $s(p) = s(e_1)$, $r(p) = r(e_n)$, and regard each vertex $v$ as a path of length zero with $s(v) = r(v) = v$. The set of paths is $\Path(E)$, and $p^* = e_n^* \cdots e_1^*$ is the \textit{ghost path} of $p$. A path $p$ with $|p| \geq 1$ and $r(p) = s(p)$ is a \textit{closed path} based at $s(p)$; it is \textit{simple} if it does not pass through its base more than once, and it is a \textit{cycle} if it does not pass through any vertex twice. An \textit{exit} of a path $e_1 \cdots e_n$ is an edge $f$ with $s(f) = s(e_i)$ and $f \neq e_i$ for some $i$. The graph $E$ satisfies Condition (L) if every cycle has an exit, and Condition (K) if no vertex is the base of exactly one simple closed path. A graph is \textit{acyclic} if it has no cycles.

\begin{definition}
Given an arbitrary graph $E$, the \textit{Leavitt path algebra
}$L_{R}(E)$ is the $R$-algebra generated by a set $\{v:v\in
E^{0}\}$ of pairwise orthogonal idempotents and a set
$\{e,e^{\ast}:e\in E^{1}\}$ satisfying:
\begin{enumerate}
\item[(1)] $s(e)e=e=er(e)$ for all $e\in E^{1}$;
\item[(2)] $r(e)e^{\ast}=e^{\ast}=e^{\ast}s(e)$\ for all $e\in E^{1}$;
\item[(3)] (CK1) $e^{\ast}e=r(e)$ and $e^{\ast}f=0$ for all $e \neq f$ in $E^1$;
\item[(4)] (CK2) $v=\sum_{e \in s^{-1}(v)}ee^{\ast}$ for every regular vertex $v$.
\end{enumerate}
\end{definition}

A family $\{v, e, e^*\}$ of elements of an $R$-algebra $B$ with the same properties is a \textit{Leavitt $E$-family} in $B$, and $L_R(E)$ is universal for them: every Leavitt $E$-family in $B$ induces a unique $R$-algebra homomorphism $L_R(E) \rightarrow B$. The algebra $L_R(E)$ has an involution $*$ extending $e \mapsto e^*$ and fixing the vertices, and it is $\Z$-graded by $\deg(v) = 0$, $\deg(e) = 1$ and $\deg(e^*) = -1$; an ideal is \textit{graded} if it contains the homogeneous components of its elements, as happens for ideals generated by homogeneous elements. The algebra $L_R(E)$ is spanned by the elements $\alpha\beta^*$ with $\alpha, \beta \in \Path(E)$ and $r(\alpha) = r(\beta)$, and, if $E^0$ is finite, it is unital with $1 = \sum_{v \in E^0} v$. We will use the following standard facts.

\begin{enumerate}
\item[(F1)] For an ideal $A$ of $L_R(E)$ and $k \in R$, the set $\{ w \in E^0 \ : \ kw \in A\}$ is hereditary and saturated (see Section 2.2), because $k\,r(e) = e^*(kw)e$ for $e \in s^{-1}(w)$, and $kw = \sum_{e \in s^{-1}(w)} e\,(k\,r(e))\,e^*$ if $w$ is regular.
\item[(F2)] $kv \neq 0$ for every $v \in E^0$ and every $0 \neq k \in R$ \cite[Proposition 3.4]{Tomforde}.
\item[(F3)] For an ideal $J$ of $R$, $L_R(E)/JL_R(E) \cong L_{R/J}(E)$ as graded algebras, where $JL_R(E)$ is the ideal generated by the elements $jv$, $j \in J$, $v \in E^0$: the universal properties give mutually inverse homomorphisms which are the identity on generators.
\item[(F4)] Let $\Lambda$ be a non-empty set. We denote by $M_\Lambda(R)$ the ring of $\Lambda \times \Lambda$ matrices $x = (x_{\mu\nu})$ over $R$ with finitely many non-zero entries, and, for $\mu, \nu \in \Lambda$, by $E_{\mu,\nu}$ the matrix unit with $1$ in the $(\mu,\nu)$-entry and $0$ elsewhere. When $\Lambda$ has $n$ elements we write $M_n(R)$, and when $\Lambda$ is countably infinite, $M_\infty(R)$. Every ideal $N$ of $M_\Lambda(R)$ is of the form $M_\Lambda(J)$, where, for any fixed $\lambda \in \Lambda$, $J = \{ s \in R \ : \ sE_{\lambda,\lambda} \in N\}$. Indeed, if $x \in N$, the identity $E_{\lambda,\mu}\,x\,E_{\nu,\lambda} = x_{\mu\nu}E_{\lambda,\lambda}$ shows that every entry of $x$ lies in $J$; conversely, if $s \in J$, then $sE_{\mu,\nu} = E_{\mu,\lambda}(sE_{\lambda,\lambda})E_{\lambda,\nu} \in N$ for all $\mu, \nu$. Moreover, reducing the entries modulo $J$ gives $M_\Lambda(R)/M_\Lambda(J) \cong M_\Lambda(R/J)$.
\end{enumerate}

\textit{Graded Uniqueness Theorem} \cite[Theorem 5.3]{Tomforde}. If $B$ is a $\Z$-graded ring and $\pi \colon L_R(E) \rightarrow B$ is a graded ring homomorphism with $\pi(kv) \neq 0$ for every $v \in E^0$ and $0 \neq k \in R$, then $\pi$ is injective.

\textit{Cuntz--Krieger Uniqueness Theorem} \cite[Theorem 2.2.16]{AAS}. If $K$ is a field, $E$ satisfies Condition (L), and $\pi \colon L_K(E) \rightarrow B$ is a ring homomorphism with $\pi(v) \neq 0$ for every $v \in E^0$, then $\pi$ is injective.

\subsection{Lattices}
We write $v \geq w$ if there is a path from $v$ to $w$, and $T(v) = \{ w \ : \ v \geq w\}$. A set $H \subseteq E^0$ is \textit{hereditary} if $v \in H$ and $v \geq w$ imply $w \in H$, and \textit{saturated} if every regular vertex $v$ with $r(s^{-1}(v)) \subseteq H$ belongs to $H$. The hereditary saturated subsets of $E^0$ form a lattice $\mathcal{H}_E$, and $\overline{X}$ denotes the smallest hereditary saturated set containing $X \subseteq E^0$. For $u \in E^0$ we write $H_u = \overline{\{u\}} = \overline{T(u)}$; note that $H_v \subseteq H_u$ if and only if $v \in H_u$.

For a hereditary set $H$, a \textit{breaking vertex} of $H$ is an infinite emitter $v \notin H$ with $0 < |s^{-1}(v) \cap r^{-1}(E^0 \backslash H)| < \infty$; the set of breaking vertices is $B_H$, and, for $v \in B_H$,
$$v^H = v - \sum_{e \in s^{-1}(v), \ r(e) \notin H} ee^*.$$
By (CK1), $v^H$ and $v - v^H$ are orthogonal self-adjoint idempotents in $vL_R(E)v$, and, for $e \in s^{-1}(v)$,
\begin{equation}\label{eq:vH}
v^He = e \ \text{ if } r(e) \in H, \qquad v^He = 0 \ \text{ if } r(e) \notin H.
\end{equation}
An \textit{admissible pair} $(H,S)$ consists of $H \in \mathcal{H}_E$ and $S \subseteq B_H$. The set $\mathcal{T}_E$ of admissible pairs is a lattice for the order $(H_1,S_1) \leq (H_2,S_2)$ if $H_1 \subseteq H_2$ and $S_1 \subseteq H_2 \cup S_2$ \cite[Proposition 2.5.6]{AAS}, in which every family has a supremum \cite[Proposition 6.3]{Rigby}. We write $\mathcal{T}_E^* = \mathcal{T}_E \backslash \{(\emptyset,\emptyset)\}$; since $B_\emptyset = \emptyset$, we have $H \neq \emptyset$ for every $(H,S) \in \mathcal{T}_E^*$.

\textbf{The set $\Ext$.} The element $v^H$ only depends on $v$ and on the set of edges of $v$ with range outside $H$, so different pairs may give the same element (see Example \ref{ex:pairs}). For this reason we consider the disjoint union
$$\Ext= E^0 \sqcup \{v^H \ : \ H \in \mathcal{H}_E, \ v \in B_H \},$$
where each $v^H$ is a formal symbol indexed by the pair $(H,v)$, which determines the element $v^H$ of $L_R(E)$. If $E$ is row-finite, then $\Ext = E^0$.

For a ring $R$, $\mathcal{L}(R)$ denotes the lattice of ideals of $R$ and, if $R$ is graded, $\mathcal{L}_{gr}(R)$ the lattice of graded ideals.

\subsection{Basic ideals and graded ideals}
For $(H,S) \in \mathcal{T}_E$, $I(H,S)$ denotes the ideal generated by $H$ and the elements $v^H$, $v \in S$ \cite[Lemma 2.4.6]{AAS}, \cite[Lemma 3.3]{Larki}; it is a graded ideal. Over a field, $(H,S) \mapsto I(H,S)$ is a lattice isomorphism onto the graded ideals \cite[Theorem 5.7]{T}. Over a commutative ring, the following notion was introduced by M. Tomforde \cite{Tomforde} for row-finite graphs and by H. Larki \cite{Larki} in general.

\begin{definition}\label{Def:basicgradedideal}
An ideal $A$ of $L_R(E)$ is \textit{basic} if, for every $0 \neq k \in R$ and every $x \in L_R(E)$ which is a vertex or an element $v^H$, $kx \in A$ implies $x \in A$. Otherwise, $A$ is \textit{non-basic}.
\end{definition}

Over a field all ideals are basic. By \cite[Theorem 3.10]{Larki}, $(H,S) \mapsto I(H,S)$ is a bijection from $\mathcal{T}_E$ onto the set of basic graded ideals of $L_R(E)$, with inverse $A \mapsto (A \cap E^0, \{ v \in B_{A \cap E^0} \ : \ v^{A \cap E^0} \in A\})$; it is a lattice isomorphism, since it and its inverse preserve the order (see (\ref{eq:vHH}) and Proposition \ref{prop:maxbasic}).

\begin{definition} \cite[Definition 6.5]{Rigby}
    A \textit{saturated function} is a map $f : \mathcal{T}_E^* \rightarrow \mathcal{L}(R)$ such that $f ( \bigvee_{i} (H_i, S_i) ) = \bigcap_i f(H_i, S_i)$ for every non-empty family $\{ (H_i, S_i) \}_i$ in $\mathcal{T}_E^*$. They form a set $\mathcal{F}_{E,R}$, ordered by $f_1 \leq f_2$ if $f_1(H,S) \subseteq f_2(H,S)$ for all $(H,S)$.
\end{definition}
Applying the definition to a pair of comparable elements shows that saturated functions are order-reversing.

\begin{theorem}\label{cor:sattogr} \cite[Theorem 6.23 and Corollary 6.27]{Rigby}
    The map $A \mapsto f_A$ is a lattice isomorphism from $\mathcal{L}_{gr}(L_R(E))$ onto $\mathcal{F}_{E,R}$, where
    $$f_A(H,S) = \{ k \in R \ : \ kv \in A \text{ for all } v \in H, \text{ and } kw^H \in A \text{ for all } w \in S \}.$$
    Its inverse sends $f$ to the ideal generated by the elements $kv$ and $kw^H$ with $k \in f(H,S)$, $v \in H$, $w \in S$ and $(H,S) \in \mathcal{T}_E^*$.
\end{theorem}

Since $f_A(H,S)$ is an ideal of $R$, it equals $R$ exactly when it contains $1$; thus
\begin{equation}\label{eq:fR}
f_A(H,S) = R \iff I(H,S) \subseteq A.
\end{equation}
Moreover, $A$ is basic if and only if $f_A$ takes values in $\{0,R\}$ \cite[Remark 6.29]{Rigby}.

\begin{example} \label{Larki}
Let $R = \Z$ and let $E$ be the graph
\[
\xymatrix{
\bullet_{u}\ar@(ul,dl)_{c} \ar@/_.3pc/[r]_{\infty} &
\bullet_{v} }
\]
\cite[Example 3.4]{Larki}. The only non-trivial hereditary saturated set is $H=\{v\}$, with $B_H=\{u \}$ and $u^H= u - cc^*$, so $\mathcal{T}^*_E= \{ (H,\emptyset), (H,B_H), (E^0, \emptyset)\}$ is a chain. On a chain, the saturated functions are the order-reversing functions. Let $A$ be the graded ideal corresponding to $f(H,\emptyset) = \Z$, $f(H,B_H) = 2\Z$ and $f(E^0,\emptyset) = 0$; it is generated by $v$ and $2u^H$. Since $1 \notin f_A(H,B_H) = 2\Z$, we have $u^H \notin A$, while $2u^H \in A$; so $A$ is non-basic.
\end{example}

The next result is the analogue for basic graded ideals of \cite[Lemma 3.5]{Ranga1}, which states that every ideal of a Leavitt path algebra over a field contains a largest graded ideal.

\begin{proposition}\label{prop:maxbasic}
        Let $E$ be an arbitrary graph, $A$ a graded ideal of $L_R(E)$, $H = A \cap E^0$ and $S = \{v \in B_H \ : \ v^H \in A\}$. Then $(H,S)$ is the maximum of the set of admissible pairs $(H',S')$ such that $H' \subseteq A$ and $v^{H'}\in A$ for all $v \in S'$, and $I(H,S)$ is the largest basic graded ideal contained in $A$.
\end{proposition}
\begin{proof}
    By (F1), $H$ is hereditary and saturated, so $(H,S)$ belongs to the set. Let $(H',S')$ be another element; then $H' \subseteq H$. Let $v \in S' \backslash H$. The edges of $v$ with range outside $H$ are among the finitely many with range outside $H'$, and
    \begin{equation}\label{eq:vHH}
    v^{H'} = v - \sum_{e \in s^{-1}(v), \ r(e) \notin H} ee^* - \sum_{e \in s^{-1}(v), \ r(e) \in H \backslash H'} e e^*.
    \end{equation}
    The last sum lies in $A$, since $ee^* = e\,r(e)\,e^*$, so $v - \sum_{r(e) \notin H} ee^* \in A$. If $v$ had no edge with range outside $H$, this element would be $v$, contradicting $v \notin A$. Hence $v \in B_H$, the element is $v^H \in A$, and $v \in S$. Thus $(H',S') \leq (H,S)$. The last claim follows: the basic graded ideals are the ideals $I(H',S')$, and if $I(H',S') \subseteq A$, then $(H',S')$ belongs to the set, so $(H',S') \leq (H,S)$ and $I(H',S') \subseteq I(H,S)$.
\end{proof}

\subsection{Arbitrary ideals and Condition (K)}
The non-graded ideals of $L_K(E)$ over a field were described by K.M. Rangaswamy \cite{Ranga}, and all ideals of $L_R(E)$ by S. Rigby and T. van den Hove \cite{Rigby}.

\begin{definition}\label{def:Du}\cite[Notation 6.12 and Definition 6.13]{Rigby}
    A cycle $C$ is a \textit{cycle without (K)} if $C$ is the only simple closed path based at $s(C)$. Two cycles are equivalent if they have the same edges in the same cyclic order, and $C_u(E)$ is the set of equivalence classes of cycles without (K). For a cycle $C$, $C^0$ is its set of vertices, and $C_{\downarrow}$ is the hereditary saturated closure of the ranges of its exits. Finally, $\mathcal{D}_{E,R}$ is the set of pairs $(f,g)$ with $f \in \mathcal{F}_{E,R}$ and $g: C_u(E) \rightarrow \mathcal{L}(R[x,x^{-1}])$ such that $g(C) \cap R = f(\overline{C^0}, \emptyset)$ and, if $C_{\downarrow} \neq \emptyset$, $g(C) \subseteq f(C_{\downarrow}, \emptyset)[x,x^{-1}]$.
\end{definition}

Note that $E$ satisfies Condition (K) if and only if $C_u(E) = \emptyset$: if a vertex is the base of a unique simple closed path $p$, then $p$ is a cycle, since a vertex repeated in $p$ would allow us to remove a closed subpath and obtain a second one.

\begin{theorem}\cite[Theorems 6.23 and 6.26]{Rigby}\label{th:Rigby626}
    The ideals of $L_R(E)$ correspond bijectively to the pairs $(f,g) \in \mathcal{D}_{E,R}$, and the ideal corresponding to $(f,g)$ is graded if and only if $g(C) = f(\overline{C^0}, \emptyset)[x,x^{-1}]$ for every $C \in C_u(E)$, or, equivalently, if and only if each $g(C)$ is a graded ideal of $R[x,x^{-1}]$.
\end{theorem}

The following result combines \cite[Theorem 7.17]{Tomforde}, \cite[Theorem 3.18]{Larki} and Theorem \ref{th:Rigby626}.

\begin{theorem}\label{th:K}
    For an arbitrary graph $E$, the following are equivalent: (a) $E$ satisfies Condition (K); (b) all ideals of $L_R(E)$ are graded; (c) all basic ideals of $L_R(E)$ are graded.
\end{theorem}
\begin{proof}
    If (K) holds, then $C_u(E) = \emptyset$ and every ideal is graded by Theorem \ref{th:Rigby626}. Clearly (b) implies (c), and (c) implies (a) by \cite[Theorem 3.18]{Larki}.
\end{proof}

\section{Graded ideal functions}\label{Sec:gradedidealfunction}
In her Master's thesis \cite{Phoebe}, P. McDougall describes the ideals of $L_R(E)$, for graphs satisfying Condition (K), by means of maps $\tau: E^0 \rightarrow \mathcal{L}(R)$. We introduce a variant for graded ideals, defined on $\Ext$, which also takes into account the breaking vertices of $E$ and describes the graded ideals through their generators.

For $\varphi : \Ext \rightarrow \mathcal{L}(R)$ and $(H,S) \in \mathcal{T}_E^*$, we write
$$
\varphi[H,S] := \bigcap_{u \in H} \varphi(u) \cap \bigcap_{v \in S} \varphi(v^H),
$$
where an empty intersection is $R$.

\begin{definition}\label{Def:varphi}
A function $\varphi: \Ext \rightarrow \mathcal{L}(R)$ is a \textit{graded ideal function of $L_R(E)$} if:
\begin{enumerate}
    \item[(GF1)] $\varphi[\bigvee_i (H_i,S_i)] = \bigcap_i \varphi[H_i,S_i]$ for every non-empty family $\{(H_i,S_i)\}_i$ in $\mathcal{T}_E^*$;
    \item[(GF2)] $\varphi(u) \subseteq \varphi(v)$ for all $u, v \in E^0$ with $H_v \subseteq H_u$;
    \item[(GF3)] $\varphi(v^H) \subseteq \bigcap_{w \in H} \varphi(w)$ for all $H \in \mathcal{H}_E$ and $v \in B_H$.
\end{enumerate}
The set $\mathcal{G}_{E,R}$ of graded ideal functions is ordered by $\varphi_1 \leq \varphi_2$ if $\varphi_1(x) \subseteq \varphi_2(x)$ for all $x \in \Ext$; its minimum and maximum are the constant functions $0$ and $R$.
\end{definition}

Since $H_w \subseteq H_u$ if and only if $w \in H_u$, conditions (GF2) and (GF3) say that $\varphi(u) = \varphi[H_u,\emptyset]$ for every $u \in E^0$ and $\varphi(v^H) = \varphi[H,\{v\}]$ for every $v^H \in \Ext$, respectively.

For $f \in \mathcal{F}_{E,R}$, define $\varphi_f : \Ext \rightarrow \mathcal{L}(R)$ by $\varphi_f(u) = f(H_u,\emptyset)$ and $\varphi_f(v^H) = f( H, \{v\})$, and, for a graded ideal $A$, write $\varphi_A = \varphi_{f_A}$, or just $\varphi$. Conversely, for $\varphi : \Ext \rightarrow \mathcal{L}(R)$, define $f_\varphi(H,S) = \varphi[H,S]$. The following example shows that $\varphi_f(v^H)$ depends on the pair $(H,v)$, and not only on the element $v^H \in L_R(E)$.

\begin{example}\label{ex:pairs}
Let $R = \Z$ and let $E$ be the graph
\[
\xymatrix{ \bullet_{w} & \bullet_{v} \ar[l]_{\infty} \ar[r]^{g} & \bullet_{x} & \bullet_{z} }
\]
where $v$ emits infinitely many edges to $w$ and one edge $g$ to $x$. Since $E$ has no regular vertices, every hereditary set is saturated. For $H_1 = \{w\}$ and $H_2 = \{w,z\}$ we have $v \in B_{H_1} \cap B_{H_2}$ and $v^{H_1} = v - gg^* = v^{H_2}$. Let $A = I(\{v,w,x\},\emptyset)$, the ideal generated by $v$. Then $\varphi_A(v^{H_1}) = \Z$, since $w, v^{H_1} \in A$, whereas $\varphi_A(v^{H_2}) = 0$, since $A$ is basic and $z \notin A$.
\end{example}

\begin{theorem} \label{th:iso}
   Let $E$ be an arbitrary graph.
   \begin{enumerate}
       \item The maps $f \mapsto \varphi_f$ and $\varphi \mapsto f_\varphi$ are mutually inverse lattice isomorphisms between $\mathcal{F}_{E,R}$ and $\mathcal{G}_{E,R}$. Consequently, $\mathcal{G}_{E,R}$ is a lattice isomorphic to $\mathcal{L}_{gr}(L_R(E))$.
       \item The graded ideal corresponding to $\varphi \in \mathcal{G}_{E,R}$ is generated by the elements $ku$ and $hv^H$ with $u \in E^0$, $k \in \varphi(u)$, $v^H \in \Ext$ and $h \in \varphi(v^H)$. Conversely, for a graded ideal $A$,
   $$
   \varphi_A(u) = \{ k \in R \ : \ ku \in A \}, \qquad \varphi_A(v^H) = \{ k \in R \ : \ k v^H \in A \mbox{ and } kH \subseteq A\}.
   $$
   \end{enumerate}
\end{theorem}
\begin{proof}
By (\ref{eq:join}) below, a saturated function is determined by its values on the pairs $(H_u,\emptyset)$ and $(H,\{v\})$, which are the values recorded by $\varphi_f$. Every $(H,S) \in \mathcal{T}_E^*$ satisfies
\begin{equation}\label{eq:join}
(H,S) = \bigvee_{u \in H} (H_u, \emptyset) \vee \bigvee_{v \in S} (H,\{v\}).
\end{equation}
Every pair on the right-hand side lies below $(H,S)$, because $H_u \subseteq H$ for $u \in H$ and $v \in S$. Conversely, if $(H',S')$ lies above all of them, then $H'$ contains every $H_u$ with $u \in H$, hence contains $H$, and every $v \in S$ belongs to $H' \cup S'$; that is, $(H,S) \leq (H',S')$. So $(H,S)$ is the least upper bound, as claimed. If $f$ is a saturated function, applying $f$ to (\ref{eq:join}) turns the supremum into an intersection, and we obtain
\begin{equation}\label{eq:phiF}
f(H,S) = \bigcap_{u \in H} f(H_u,\emptyset) \cap \bigcap_{v \in S} f(H,\{v\}) = \varphi_f[H,S] \qquad \text{for every } (H,S) \in \mathcal{T}_E^*.
\end{equation}

Now we prove (1). Identity (\ref{eq:phiF}) says that $f_{\varphi_f} = f$, and it also shows that $\varphi_f$ is a graded ideal function. Indeed, through (\ref{eq:phiF}), condition (GF1) for $\varphi_f$ is exactly the saturation of $f$. Condition (GF2) holds because $H_v \subseteq H_u$ gives $(H_v,\emptyset) \leq (H_u,\emptyset)$, and $f$ is order-reversing. Condition (GF3) is (\ref{eq:phiF}) for the pair $(H,\{v\})$, which reads $\varphi_f(v^H) = \bigcap_{w \in H}\varphi_f(w) \cap \varphi_f(v^H)$. In the other direction, let $\varphi$ be a graded ideal function. Condition (GF1) says precisely that $f_\varphi$ is saturated, and, as observed after Definition \ref{Def:varphi}, conditions (GF2) and (GF3) say that $\varphi(u) = f_\varphi(H_u,\emptyset)$ and $\varphi(v^H) = f_\varphi(H,\{v\})$, that is, $\varphi_{f_\varphi} = \varphi$. Hence the two maps are mutually inverse bijections. Both of them preserve the order, since they only take values and intersections of values, so they are lattice isomorphisms. Composing with the isomorphism of Theorem \ref{cor:sattogr}, we obtain $\mathcal{G}_{E,R} \cong \mathcal{L}_{gr}(L_R(E))$.

For (2), let $A$ be the graded ideal corresponding to $\varphi$, so that $f_A = f_\varphi$. By Theorem \ref{cor:sattogr}, $A$ is generated by the elements $ku$ and $kv^H$ with $(H,S) \in \mathcal{T}_E^*$, $k \in f_\varphi(H,S)$, $u \in H$ and $v \in S$. We compare this set of generators with the one in the statement. On the one hand, the generators $ku$ with $k \in \varphi(u)$ and $hv^H$ with $h \in \varphi(v^H)$ are among the former, as they come from the pairs $(H_u,\emptyset)$ and $(H,\{v\})$. On the other hand, if $k \in f_\varphi(H,S)$, $u \in H$ and $v \in S$, then $k \in \varphi(u)$ and $k \in \varphi(v^H)$, because $f_\varphi$ is order-reversing and $(H_u,\emptyset), (H,\{v\}) \leq (H,S)$. So both sets generate the same ideal. Finally, by definition, $\varphi_A(u) = f_A(H_u,\emptyset)$ is the set of $k$ with $kw \in A$ for every $w \in H_u$. This set is contained in $\{k \ : \ ku \in A\}$, and the reverse inclusion holds because, if $ku \in A$, the set $\{w \ : \ kw \in A\}$ is hereditary and saturated by (F1) and contains $u$, so it contains $H_u$. The formula for $\varphi_A(v^H) = f_A(H,\{v\})$ is the one given in Theorem \ref{cor:sattogr}.
\end{proof}

\begin{corollary}\label{Lem:varphi2}
    Let $\varphi$ be a graded ideal function.
    \begin{enumerate}
        \item $\varphi(u) = \bigcap_{w \in H_u} \varphi(w)$ for every $u \in E^0$.
        \item $\varphi(v) \cap \bigcap_{u \in H} \varphi(u) \subseteq \varphi(v^H)$ for every $v^H \in \Ext$.
        \item If $x \geq y$ in $E^0$, then $\varphi(x) \subseteq \varphi(y)$; and, by (GF2), if $H_u = H_v$, then $\varphi(u) = \varphi(v)$.
        \item If $x \geq y \geq z$ and $\varphi(x) = \varphi(z)$, then $\varphi(y)=\varphi(x)$.
    \end{enumerate}
\end{corollary}
\begin{proof}
(1) is (GF2). For (2), write $\varphi = \varphi_A$ and use Theorem \ref{th:iso}(2): if $kv \in A$ and $kH \subseteq A$, then $kv^H = (kv)v^H \in A$. For (3), $x \geq y$ gives $y \in H_x$, so $H_y \subseteq H_x$, and we apply (GF2); the second statement is (GF2) applied twice. Finally, (4) follows from (3).
\end{proof}

Condition (GF1) involves suprema of arbitrary families in $\mathcal{T}_E^*$. For row-finite graphs, it can be replaced by a local condition, in the spirit of \cite{Phoebe}.

\begin{proposition}\label{prop:local}
Let $E$ be a row-finite graph, so that $\Ext = E^0$. A function $\varphi: E^0 \rightarrow \mathcal{L}(R)$ is a graded ideal function if and only if
\begin{equation}\label{eq:local}
\varphi(v) = \bigcap_{e \in s^{-1}(v)} \varphi(r(e)) \quad \text{ for every regular vertex } v.
\end{equation}
\end{proposition}
\begin{proof}
If $\varphi = \varphi_A$, then (\ref{eq:local}) follows from Theorem \ref{th:iso}(2), since, for a regular vertex $v$, $kv \in A$ if and only if $k\,r(e) \in A$ for every $e \in s^{-1}(v)$, by (F1). Conversely, assume (\ref{eq:local}) and, for $k \in R$, let $H_k = \{ u \in E^0 \ : \ k \in \varphi(u)\}$. By (\ref{eq:local}), $\varphi(s(e)) \subseteq \varphi(r(e))$ for every edge $e$, so $H_k$ is hereditary, and it is saturated, again by (\ref{eq:local}). Since $E$ is row-finite, $\mathcal{T}_E^*$ consists of the pairs $(H,\emptyset)$, and $\bigvee_i (H_i,\emptyset) = (\overline{\bigcup_i H_i},\emptyset)$. As $k \in \varphi[H,\emptyset]$ means $H \subseteq H_k$, and $H_k$ is hereditary and saturated, $k$ belongs to either side of (GF1) if and only if $\bigcup_i H_i \subseteq H_k$. Similarly, if $k \in \varphi(u)$, then $H_u \subseteq H_k$, which gives (GF2); and (GF3) is empty.
\end{proof}

\begin{example} \label{ex:ToeplitzGF}
We will use the Toeplitz graph
\[ T:
\xymatrix{
\bullet_{u}\ar@(ul,dl)_{c} \ar@/_.1pc/[r]_e &
\bullet_{v}}
\]
as a running example. Its only non-trivial hereditary saturated set is $\{v\}$, so $\mathcal{T}_{T}^*$ consists of $(\{v\},\emptyset)$ and $(T^0,\emptyset)$, with $H_v = \{v\}$ and $H_u = T^0$. By Proposition \ref{prop:local}, the graded ideal functions are given by $\varphi(u) = \varphi(u) \cap \varphi(v) \subseteq \varphi(v)$. Hence the graded ideals of $L_R(T)$ correspond to the pairs of ideals $I \subseteq J$ of $R$, via $J = \varphi(v)$ and $I = \varphi(u)$, and the graded ideal corresponding to $(I,J)$ is generated by $Jv$ and $Iu$. For $R = \Z$, we write $J = a\Z$ and $I = b\Z$, where $a$ divides $b$.
\end{example}

\section{Quotients by \texorpdfstring{$I$}{I}-basic graded ideals}\label{Sec:Main_via_f}

The \textit{quotient graph} of $E$ by an admissible pair $(H,S)$ is the graph $E\backslash(H,S)$ with
$$(E\backslash(H,S))^{0}=(E^{0}\backslash H)\cup\{v^{\prime}:v\in B_{H}\backslash S\},$$
$$(E\backslash(H,S))^{1}=\{e\in E^{1}:r(e)\notin H\}\cup\{e^{\prime}:r(e)\in B_{H}\backslash S\},$$
where $s(e^{\prime})=s(e)$ and $r(e^{\prime})=r(e)^{\prime}$. It is known that $L_{R}(E)/I(H,S) \cong L_{R}(E\backslash (H,S))$ \cite[Theorem 5.7]{T}, \cite[Theorem 3.10]{Larki}. If $A$ is a graded ideal, $H = A \cap E^0$ and $S = \{v \in B_H \ : \ v^H \in A\}$, then $I(H,S) \subseteq A$ by Proposition \ref{prop:maxbasic}, so $L_R(E)/A$ is a quotient of $L_R(E\backslash(H,S))$. The next proposition refines this.

\begin{proposition}\label{Th:Isom}
Let $E$ be an arbitrary graph, $A$ a graded ideal of $L_R(E)$ with saturated function $f$, $H=A \cap E^0$, $S= \{v \in B_H \ : \ v^H \in A \}$ and $I=f(E^0,\emptyset)$. Then $L_R(E)/A$ is an $R/I$-algebra, and there is a graded $R/I$-algebra epimorphism
        $$\psi : L_{R/I}(E\backslash(H,S)) \rightarrow L_{R}(E)/A, \qquad \psi(x) = P_x + A \ \text{ for every vertex } x,$$
where $P_u = u$ if $u \in E^0 \backslash (H \cup B_H)$ or $u \in S$, and $P_u = u - u^H$, $P_{u'} = u^H$ if $u \in B_H \backslash S$.
\end{proposition}

\begin{proof}
Write $L = L_R(E)$ and $F = E\backslash(H,S)$. We first check that $IL \subseteq A$, so that $L/A$ is an $R/I$-algebra. If $h \in I = f(E^0,\emptyset)$, then $hw \in A$ for every vertex $w$, by the definition of $f$. Since $L$ is spanned by the elements $\alpha\beta^*$ and $h\alpha\beta^* = \alpha(h\,r(\alpha))\beta^*$, we get $hL \subseteq A$.

By the universal property of $L_{R/I}(F)$, it suffices to find a Leavitt $F$-family in $L/A$. Let $u \in E^0 \backslash H$. If $u \in B_H \backslash S$, then $u = u^H + (u - u^H)$ is a sum of orthogonal idempotents, with $u^He = e$ if $r(e) \in H$ and $(u-u^H)e = e$ if $r(e) \notin H$; in $F$ these correspond to the new sink $u'$ and to the vertex $u$, which is regular in $F$. This explains the choice of the elements $P_x$ in the statement. For the edges, let $e \in E^1$ with $r(e) \notin H$, and let $g$ be $e$, or $e'$ when $r(e) \in B_H \backslash S$; we set $D_g = eP_{r(g)}$. Since $A$ is generated by self-adjoint elements (Theorem \ref{th:iso}(2)), we have $A^* = A$, and $D_g^*$ will play the role of $g^*$.

We check the relations modulo $A$. The elements $P_x$ are pairwise orthogonal idempotents, because $u^H$ and $u - u^H$ are. Relations (1) and (2) reduce to $P_{s(e)}e = e$, which is clear when $P_{s(e)} = s(e)$; when $s(e) \in B_H \backslash S$, the edge $e$ is one of the edges in the sum defining $P_{s(e)}$, so $P_{s(e)}e = e$ by (CK1). For (CK1), let $g$ and $h$ be edges of $F$ with underlying edges $e$ and $f$. Then $D_g^*D_h = P_{r(g)}e^*fP_{r(h)}$, which is $0$ if $e \neq f$; if $e = f$, it equals $P_{r(g)}\,r(e)\,P_{r(h)} = P_{r(g)}P_{r(h)}$, which is $P_{r(g)}$ if $g = h$ and $0$ otherwise. It remains to verify (CK2). Note that $F$ has new sinks and that some infinite emitters of $E$ become regular in $F$. A vertex $u$ of $F$ is regular exactly when it is a vertex of $E\backslash H$ emitting finitely many, and at least one, edges with range outside $H$. Its edges in $F$ are these edges $e$, together with the edges $e'$ for $r(e) \in B_H \backslash S$, and since $P_{r(e)} + P_{r(e)'} = r(e)$ in the latter case,
$$\sum_{g \in s_F^{-1}(u)} D_gD_g^* = \sum_{e \in s^{-1}(u), \, r(e) \notin H} ee^*.$$
It remains to compare this element with $P_u$ modulo $A$. If $u$ is regular in $E$, then $u = \sum_{e \in s^{-1}(u)} ee^*$ by (CK2), and the terms with $r(e) \in H$ belong to $A$, because $ee^* = e\,r(e)\,e^*$; so our element is congruent to $u = P_u$. If $u$ is not regular in $E$, it is an infinite emitter with finitely many edges leaving $H$, that is, $u \in B_H$, and our element is $u - u^H$. This is $P_u$ when $u \notin S$, and it is congruent to $u = P_u$ when $u \in S$, because then $u^H \in A$.

The universal property now gives an $R/I$-algebra homomorphism $\psi$ with $\psi(x) = P_x + A$, $\psi(g) = D_g + A$ and $\psi(g^*) = D_g^* + A$. It is graded, since the elements $P_x$ have degree $0$ and the elements $D_g$ have degree $1$. Finally, $\psi$ is onto, because every generator of $L/A$ lies in its image: the vertices of $H$ and the edges with range in $H$ belong to $A$; a vertex $u \notin H \cup (B_H\backslash S)$ is $P_u$, and a vertex $u \in B_H\backslash S$ is $P_u + P_{u'}$; an edge $e$ with $r(e) \notin H$ is $D_e$, or $D_e + D_{e'}$ when $r(e) \in B_H\backslash S$; and the ghost edges are handled by taking adjoints.
\end{proof}

\begin{definition}\label{Def:I-basicgradedideal} Let $I$ be an ideal of $R$. A graded ideal $A$ of $L_R(E)$ with saturated function $f$ is \textit{$I$-basic} if $\im(f) = \{ I, R\}$ or $\im(f) = \{ I\}$.
\end{definition}

By (\ref{eq:phiF}), every value of $f$ is an intersection of values of $\varphi$, and every value of $\varphi$ is a value of $f$. Hence $A$ is $I$-basic if and only if $\varphi$ only takes the values $I$ and $R$.

Since $(E^0,\emptyset)$ is the maximum of $\mathcal{T}_E^*$ and $f$ is order-reversing, $f(E^0,\emptyset)$ is the smallest value of $f$; so, if $A$ is $I$-basic, then $I = f(E^0,\emptyset)$.

\begin{remark}\label{Rem:Basic} The $I$-basic graded ideals generalize the basic ones: a graded ideal $A \neq L_R(E)$ is basic if and only if it is $0$-basic, because some vertex is not in $A$, so $f_A(E^0,\emptyset) \neq R$. If $I$ is a non-zero proper ideal, every $I$-basic graded ideal $A$ is non-basic: there is a vertex $w \notin A$, while $Iw \subseteq A$.
\end{remark}

The $I$-basic graded ideals admit an explicit description.

\begin{lemma}\label{lem:Ibasic}
Let $I \neq R$ be an ideal of $R$ and $(H,S) \in \mathcal{T}_E$ with $(H,S) \neq (E^0,\emptyset)$. Then $A = I(H,S) + IL_R(E)$ is an $I$-basic graded ideal, $A \cap E^0 = H$, $\{v \in B_H \ : \ v^H \in A\} = S$, and $f_A(H',S') = R$ if $(H',S') \leq (H,S)$, while $f_A(H',S') = I$ otherwise. Conversely, every $I$-basic graded ideal $A \neq L_R(E)$ is of this form, with $H = A \cap E^0$ and $S = \{v \in B_H \ : \ v^H \in A\}$.
\end{lemma}
\begin{proof}
We reduce modulo $I$, where $A$ becomes a basic graded ideal. By (F3), the quotient map $\pi \colon L_R(E) \rightarrow L_{R/I}(E)$ is onto, with kernel $IL_R(E)$, so the ideals of $L_R(E)$ containing $IL_R(E)$ correspond to the ideals of $L_{R/I}(E)$. Under this correspondence, $A = I(H,S) + IL_R(E)$ corresponds to the ideal of $L_{R/I}(E)$ generated by $H$ and the elements $v^H$, $v \in S$, which is the basic graded ideal $I' = I(H,S)$ of $L_{R/I}(E)$. In particular, $A = \pi^{-1}(I')$, and $A$ is graded.

To compute $f_A$, we need to know when $kx \in A$, for $x$ a vertex or an element $w^{H'}$. If $k \in I$, this always happens. If $k \notin I$, then $kx \in A$ if and only if $(k+I)\pi(x) \in I'$, and, since $I'$ is basic and $k + I \neq 0$, this happens if and only if $\pi(x) \in I'$, that is, if and only if $x \in A$. Consequently, $f_A(H',S')$ is $R$ if the generators of $I(H',S')$ belong to $A$, and $I$ otherwise. By \cite[Theorem 3.10]{Larki} over $R/I$, the generators of $I(H',S')$ belong to $A$ if and only if $I(H',S') \subseteq I'$ in $L_{R/I}(E)$, if and only if $(H',S') \leq (H,S)$. This gives the formula for $f_A$. Since $(H,S) \neq (E^0,\emptyset)$, the value $I$ is attained at $(E^0,\emptyset)$, and $A$ is $I$-basic. The same theorem gives $I' \cap E^0 = H$ and $\{v \in B_H \ : \ v^H \in I'\} = S$, and, as we have just seen with $k = 1$, a vertex or an element $v^H$ lies in $A$ if and only if its image lies in $I'$; so $A \cap E^0 = H$ and $\{v \in B_H \ : \ v^H \in A\} = S$.

Conversely, let $A \neq L_R(E)$ be an $I$-basic graded ideal, and let $H = A \cap E^0$ and $S = \{v \in B_H \ : \ v^H \in A\}$. Then $I \neq R$, and $(H,S) \neq (E^0,\emptyset)$, since otherwise $A$ would contain all the vertices. We show that $f_A$ is the function computed above for $B = I(H,S) + IL_R(E)$; since graded ideals are determined by their saturated functions (Theorem \ref{cor:sattogr}), this gives $A = B$. As $f_A$ only takes the values $I$ and $R$, it suffices to know where it takes the value $R$. By (\ref{eq:fR}), $f_A(H',S') = R$ if and only if $I(H',S') \subseteq A$, and, by Proposition \ref{prop:maxbasic}, this happens if and only if $(H',S') \leq (H,S)$. So $f_A = f_B$.
\end{proof}

We can now prove the main result of this section.

\begin{theorem}\label{cor:main}
Let $E$ be an arbitrary graph and $A$ a graded ideal of $L_R(E)$, with $H$, $S$, $I$ and $\psi$ as in Proposition \ref{Th:Isom}. Then $\psi$ is an isomorphism if and only if $A$ is $I$-basic. In particular, if $A$ is $I$-basic, then
 $$L_{R}(E)/A \cong L_{R/I}(E\backslash (H,S)).$$
\end{theorem}
\begin{proof}
If $I = R$, then $A = L_R(E)$ and both algebras are zero, so we assume that $I \neq R$. Then $(H,S) \neq (E^0,\emptyset)$, since otherwise $A$ would contain every vertex. We compare $A$ with $B = I(H,S) + IL_R(E)$. By Proposition \ref{prop:maxbasic}, $I(H,S) \subseteq A$, and $IL_R(E) \subseteq A$, as shown at the beginning of the proof of Proposition \ref{Th:Isom}; so $B \subseteq A$. By Lemma \ref{lem:Ibasic}, $B$ is an $I$-basic graded ideal with $B \cap E^0 = H$, $\{v \in B_H \ : \ v^H \in B\} = S$ and $f_B(E^0,\emptyset) = I$. Hence Proposition \ref{Th:Isom} also applies to $B$, with the same $H$, $S$ and $I$, and gives a graded epimorphism $\psi_B \colon L_{R/I}(E\backslash(H,S)) \rightarrow L_R(E)/B$ with $\psi_B(x) = P_x + B$. Since $\psi$ and $\psi_B$ are defined by the same elements, $\psi = q \circ \psi_B$, where $q \colon L_R(E)/B \rightarrow L_R(E)/A$ is the canonical projection.

We show that $\psi_B$ is injective. By the Graded Uniqueness Theorem over $R/I$, it suffices to see that $kP_x \notin B$ for every vertex $x$ of $E\backslash(H,S)$ and every $k \in R \backslash I$. As shown in the proof of Lemma \ref{lem:Ibasic}, if $k \notin I$ and $y$ is a vertex or an element $v^H$, then $ky \in B$ if and only if $y \in B$. This settles the vertices $x$ with $P_x = u$, $u \notin H$, and $P_x = u^H$, $u \in B_H \backslash S$, because then $P_x \notin B$. If $P_x = u - u^H$ with $u \in B_H \backslash S$, we choose $e \in s^{-1}(u)$ with $r(e) \notin H$; if $k(u - u^H)$ were in $B$, then $k\,r(e) = e^*k(u-u^H)e$ would be in $B$ by (\ref{eq:vH}), so $r(e) \in B \cap E^0 = H$, a contradiction.

Consequently, $\psi$ is an isomorphism if and only if $q$ is injective, that is, if and only if $A = B$. If $A = B$, then $A$ is $I$-basic by Lemma \ref{lem:Ibasic}. Conversely, if $A$ is $I$-basic, then $A = I(H,S) + IL_R(E) = B$, by the same lemma.
\end{proof}

For basic graded ideals, Theorem \ref{cor:main} is the quotient isomorphism of \cite[Theorem 3.10]{Larki}.

\begin{example} \label{ex:Toeplitz}
For the Toeplitz graph $T$ of Example \ref{ex:ToeplitzGF} and $R = \Z$, the graded ideal $A$ with $\varphi(v) = a\Z$ and $\varphi(u) = b\Z$ has $\im(f) = \{a\Z, b\Z\}$, and $H = \{v\}$ if $a = 1$, $H = \emptyset$ otherwise. The graph $T\backslash(\{v\},\emptyset)$ is a single loop, whose Leavitt path algebra over $R/I$ is $(R/I)[x,x^{-1}]$ (Lemma \ref{lem:structure}(4)). For $A \neq 0, L_\Z(T)$ and $n > 1$, Theorem \ref{cor:main} gives:
\begin{center}
\begin{tabular}{|c|c|c|c|}
\hline
$(a,b)$ & type of $A$ & $(H,S)$ & $L_\Z(T)/A$ \\
\hline
$(1,0)$ & $0$-basic, $A = I(\{v\},\emptyset)$ & $(\{v\},\emptyset)$ & $\Z[x,x^{-1}]$ \\
$(1,n)$ & $n\Z$-basic & $(\{v\},\emptyset)$ & $\Z_n[x,x^{-1}]$ \\
$(n,n)$ & $n\Z$-basic & $(\emptyset,\emptyset)$ & $L_{\Z_n}(T)$ \\
$a > 1$, $a \neq b$ & not $I$-basic & & see Example \ref{ex:Toeplitz2} for $(2,4)$ \\
\hline
\end{tabular}
\end{center}
\end{example}

\begin{example}\label{ex:infinityToplitz}
Let $E$ be the graph of Example \ref{Larki}:
\[ E:
\xymatrix{
\bullet_{u}\ar@(ul,dl)_{c} \ar@/_.3pc/[r]_{\infty} &
\bullet_{v} }
\]
and $R = \Z$. Since $\mathcal{T}_{E}^*$ is the chain $(\{v\}, \emptyset) \leq (\{v\}, \{u \}) \leq (E^0, \emptyset)$, the saturated functions are given by their values $a\Z \supseteq b\Z \supseteq d\Z$ on these pairs, that is, by integers $a, b, d \geq 0$ with $a$ dividing $b$ and $b$ dividing $d$. Let $A$ be the corresponding graded ideal. Then $v \in A$ if and only if $a = 1$, and $u^{\{v\}} \in A$ if and only if $b = 1$. The quotient graphs of $E$ by the first two pairs are
\[ E\backslash (\{v\}, \emptyset):
\xymatrix{
\bullet_{u}\ar@(ul,dl)_{c} \ar[r]^{c'} & \bullet_{u'}
}
 \qquad
 E\backslash (\{v\}, \{u\}):
\xymatrix{
\bullet_{u}\ar@(ul,dl)_{c}
}
\]
that is, the Toeplitz graph $T$ and a single loop, whose Leavitt path algebra over $R/I$ is $(R/I)[x,x^{-1}]$ (Lemma \ref{lem:structure}(4)). Theorem \ref{cor:main} gives, for $n > 1$:
\begin{center}
\begin{tabular}{|c|c|c|c|}
\hline
$(a,b,d)$ & type of $A$ & $(H,S)$ & $L_\Z(E)/A$ \\
\hline
$(1,0,0)$ & $0$-basic, $A = I(\{v\},\emptyset)$ & $(\{v\},\emptyset)$ & $L_\Z(T)$ \\
$(1,1,0)$ & $0$-basic, $A = I(\{v\},\{u\})$ & $(\{v\},\{u\})$ & $\Z[x,x^{-1}]$ \\
$(1,n,n)$ & $n\Z$-basic & $(\{v\},\emptyset)$ & $L_{\Z_n}(T)$ \\
$(1,1,n)$ & $n\Z$-basic & $(\{v\},\{u\})$ & $\Z_n[x,x^{-1}]$ \\
$(n,n,n)$ & $n\Z$-basic & $(\emptyset,\emptyset)$ & $L_{\Z_n}(E)$ \\
\hline
\end{tabular}
\end{center}
In all other cases with $A \neq 0$ and $A \neq L_\Z(E)$, the ideal $A$ is not $I$-basic.
\end{example}

\begin{example}\label{ex:clock}
    Consider the infinite clock graph
    \[ C_\N:
\xymatrix{
 & \bullet_{u_1}& \bullet_{u_2}\\
 &\bullet_{v} \ar[r]_{e_3} \ar[u]_{e_1} \ar[ur]_{e_2}  \ar[dl]_{\ddots} \ar[d] \ar[dr]_{e_4}&   \bullet_{u_3} \\
 & & \bullet_{u_4}}
\]
and $R = \Z$. Let $U = \{u_i \ : \ i \in \N\}$. Since $C_\N$ has no regular vertices, every $H \subseteq U$ is hereditary and saturated, with $B_H = \{v\}$ if $0 < |U \backslash H| < \infty$ and $B_H = \emptyset$ otherwise. Let $n > 1$ and $H_0 = U \backslash \{u_1\}$, so that $v^{H_0} = v - e_1e_1^*$. By Lemma \ref{lem:Ibasic}, $A_1 = I(H_0,\emptyset) + nL_\Z(C_\N)$ and $A_2 = I(H_0,\{v\}) + nL_\Z(C_\N)$ are $n\Z$-basic, with $(H,S) = (H_0,\emptyset)$ and $(H_0,\{v\})$, respectively; for instance, $f_{A_1}(H',\{v\}) = n\Z$ for every $H'$, while $f_{A_2}(H',\{v\}) = \Z$ if and only if $u_1 \notin H'$. The quotient graphs are
   \[ C_\N \backslash (H_0, \emptyset):
\xymatrix{
\bullet_{v}\ar[r]_{e_1} & \bullet_{u_1}&   \bullet_{v'}
}
\qquad
C_\N \backslash (H_0, \{v\}):
\xymatrix{
\bullet_{v}\ar[r]_{e_1} & \bullet_{u_1}
}
\]
so, by Theorem \ref{cor:main} and Lemma \ref{lem:structure} (2) and (3), $L_{\Z}(C_\N)/A_1 \cong M_2(\Z_n) \oplus \Z_n$ and $L_{\Z}(C_\N)/A_2 \cong M_2(\Z_n)$.
\end{example}

\section{Quotients for finite graphs whose cycles have no exits}\label{Sec:Main_via_varphi}
In this section, we describe the quotient of $L_R(E)$ by an arbitrary (basic or non-basic) graded ideal when $E$ is a finite graph in which no cycle has an exit. Cycles are considered up to equivalence (Definition \ref{def:Du}), and each cycle $C$ is regarded as a closed path based at a fixed vertex $v_C \in C^0$. Then $L_R(E)$ is a direct sum of matrix algebras over $R$ and $R[x,x^{-1}]$, one for each sink and each cycle (see \cite{AAPSM}, \cite[Theorems 1.8 and 3.8]{AAPSM2} for fields and \cite[Proposition 15]{KO} for commutative rings), and the quotient is a direct sum of Leavitt path algebras over quotient rings of $R$ of the graphs $_XE$ defined below. Examples \ref{ex:tree} and \ref{ex:Toeplitz2} show that neither the finiteness of $E$ nor the hypothesis on the cycles can be removed.

\subsection{The graph \texorpdfstring{$_XE$}{XE}}
A similar construction for admissible pairs, the porcupine graph \cite[Definition 3.1]{Vas}, is used in \cite{Vas} to show that every graded ideal of a Leavitt path algebra is graded isomorphic to a Leavitt path algebra.

\begin{definition} \label{Gporcupine}
Let $E$ be a row-finite graph and $X \subseteq E^0$. Let $F(X)$ be the set of paths $\alpha = e_1 \cdots e_n$, $n \geq 1$, with $r(e_n) \in X$ and $s(e_1), r(e_1), \ldots, r(e_{n-1}) \notin X$; note that if $\alpha = e\beta \in F(X)$ with $|\beta| \geq 1$, then $\beta \in F(X)$. For each $\alpha \in F(X)$, take a new vertex $w^\alpha$ and a new edge $f^\alpha$ with $s(f^\alpha)=w^\alpha$, and $r(f^\alpha) = r(\alpha)$ if $|\alpha| = 1$, $r(f^\alpha)=w^\beta$ if $\alpha = e\beta$ with $e \in E^1$ and $|\beta| \geq 1$. The graph $_XE$ has vertices $X \cup \{w^\alpha \ : \ \alpha \in F(X)\}$ and edges $\{ e \in E^1 \ : \  s(e), r(e) \in X\} \cup \{f^\alpha \ : \ \alpha \in F(X) \}$.
\end{definition}

\begin{example} \label{ex:porcurpine1}
Let $E$ be the graph
    \[ E:
\xymatrix{ & & \bullet_{v_2}\\
\bullet_{u} \ar[r]_{e_1} & \bullet_{v_1} \ar[r]_{e_3} \ar[ur]^{e_2} & \bullet_{v_3}}
\]
For $X = \{v_1\}$, $F(X) = \{ e_1 \}$ and $_XE$ is $\xymatrix{ \bullet_{w^{e_1}}\ar[r]^{f^{e_1}} & \bullet_{v_1} }$. For $X_i = \{v_i\}$, $i=2,3$, $F(X_i) = \{  e_i, e_1e_i \}$ and
\[ _{X_i}E: \xymatrix{\bullet_{w^{e_1e_i}} \ar[r]^{f^{e_1e_i}} & \bullet_{w^{e_i}}\ar[r]^{f^{e_i}} & \bullet_{v_i} }. \]
For $X = \{v_2, v_3\}$, $_XE$ is the disjoint union of $_{X_2}E$ and $_{X_3}E$.
\end{example}

\subsection{The main result}
Let $E$ be a finite graph in which no cycle has an exit; thus each vertex of a cycle emits only the edge of the cycle starting at it. Following L. Va\v{s} \cite[Definition 5.5]{VasPQ}, a vertex is \textit{terminal} if it is a sink, lies on a cycle without exits, lies on an extreme cycle, or lies on a terminal path. In our setting there are no extreme cycles, since these have exits, and no terminal paths, since every infinite path eventually goes around a cycle; so the terminal vertices are the sinks and the vertices on cycles. We denote their set by $E^0_{\mathrm{ter}}$. Va\v{s} also groups the terminal vertices into \textit{clusters} \cite[Definition 5.5]{VasPQ}, and by \cite[Lemma 5.6]{VasPQ} the cluster of a sink $w$ is $\{w\}$, and the cluster of a vertex of a cycle $C$ without exits is $C^0$. We denote by $\mathcal{Y}_E$ the set of clusters of $E$. In the language of \cite{KO, KO2}, where the sinks and the cycles are ordered by the existence of paths, the clusters correspond to the sinks and the cycles without exits, which are the minimal elements of this order. The following properties are immediate from the absence of exits:
\begin{enumerate}
\item[(T1)] a path that reaches a terminal vertex $t$ stays in the cluster of $t$; in particular, $T(t)$ is that cluster, and distinct clusters are disjoint;
\item[(T2)] non-terminal vertices are regular and lie on no closed path, since a shortest closed path based at a vertex is a cycle.
\end{enumerate}
For a sink $w$, let $P_w$ be the set of paths ending at $w$ (including $w$), and $n(w) = |P_w|$. For a cycle $C$, let $P_C$ be the set of paths $p$ ending at $v_C$ which are not of the form $qC$, and $m(C) = |P_C|$. By (T1) and (T2), a path in $P_w$ has pairwise distinct vertices, and a path in $P_C$ consists of such a path followed by less than one turn around $C$; so these sets are finite. A graded ideal function $\varphi$ is constant on each cluster $Y$, because $H_u = H_{u'}$ for $u, u'$ on the same cycle (Corollary \ref{Lem:varphi2}(3)); we write $\varphi(Y)$ for its value. Finally, for $k \in \Z$, we set $C^{\langle k \rangle} = C^k$ if $k > 0$, $C^{\langle 0 \rangle} = v_C$ and $C^{\langle k \rangle} = (C^*)^{-k}$ if $k < 0$, an element of degree $k|C|$.

The next lemma collects known results. Part (1) follows
from repeated applications of (CK2), as in the proof of
\cite[Proposition 15]{KO}. Part (2) follows from
\cite[Proposition 15]{KO}, since $L_{\{w\}}=(w)$ and
$L_{C^0}=(v_C)$ (for a field and $E$ finite, (2) is also \cite[Theorems 1.8 and 3.8]{AAPSM2}, and a special case of \cite[Proposition 6.1]{VasPQ}), and parts (3) and (4) are \cite[Proposition 14(iv) and (vi)]{KO} (see \cite[Theorem 2.6.17 and Lemma 2.7.1]{AAS} for a field).

\begin{lemma}\label{lem:structure}
Let $E$ be a finite graph in which no cycle has an exit, and let $P$ be the union of the sets $P_w$, $w$ a sink, and $P_C$, $C$ a cycle.
\begin{enumerate}
    \item $u = \sum_{p \in P, \ s(p) = u} pp^*$ for every $u \in E^0$.
    \item $L_R(E) = \bigoplus_{Y \in \mathcal{Y}_E} L_Y$, where $L_Y$ is the ideal generated by $Y$.
    \item For a sink $w$, $E_{p,q} \mapsto pq^*$ ($p,q \in P_w$) defines an isomorphism $M_{n(w)}(R) \cong L_{\{w\}}$.
    \item For a cycle $C$, $x^kE_{p,q} \mapsto pC^{\langle k \rangle}q^*$ ($p,q \in P_C$, $k \in \Z$) defines an isomorphism $$M_{m(C)}(R[x,x^{-1}]) \cong L_{C^0}.$$
\end{enumerate}
\end{lemma}

In particular, $L_{\{w\}}$ and $L_{C^0}$ are unital, with identities $\varepsilon_{\{w\}} = \sum_{p \in P_w} pp^*$ and $\varepsilon_{C^0} = \sum_{p \in P_C} pp^*$, the images of the identity matrices. Summing (1) over all vertices gives $\sum_Y \varepsilon_Y = \sum_u u = 1$, so, by (2), the elements $\varepsilon_Y$ are pairwise orthogonal central idempotents of $L_R(E)$. Moreover, since $pp^*u = pp^*$ if $s(p) = u$ and $pp^*u = 0$ otherwise,
\begin{equation}\label{eq:epsu}
\varepsilon_Y u = \sum_{p \in P_Y, \ s(p) = u} pp^* \qquad (u \in E^0,\ Y \in \mathcal{Y}_E),
\end{equation}
where $P_Y = P_w$ if $Y = \{w\}$ and $P_Y = P_C$ if $Y = C^0$. Under the isomorphisms (3) and (4), $\varepsilon_Y u$ is a diagonal matrix with entries $0$ and $1$, which is non-zero exactly when $u \geq t$ for some $t \in Y$.

\begin{theorem}\label{Th:MainI}
     Let $E$ be a finite graph in which no cycle has an exit, $A$ a graded ideal of $L_R(E)$ and $\varphi$ its graded ideal function. Then:
     \begin{enumerate}
         \item $\varphi(u) = \bigcap \{ \varphi(t) \ : \ t \in E^0_{\mathrm{ter}}, \ u \geq t\}$ for every $u \in E^0$;
         \item for $I \in \varphi(E^0_{\mathrm{ter}})$, let $X_I := \varphi^{-1}(I) \cap E^0_{\mathrm{ter}}$; then
    $$L_R(E)/A \cong \bigoplus_{I\in \varphi(E^0_{\mathrm{ter}}) \backslash \{R\}} L_{R/I}\left( _{X_I}E\right) \cong \bigoplus_{w} M_{n(w)}\left(R/\varphi(w)\right) \oplus \bigoplus_{C} M_{m(C)}\big((R/\varphi(C^0))[x,x^{-1}]\big),$$
    where $w$ runs over the sinks and $C$ over the cycles of $E$ up to equivalence;
         \item $A \mapsto (\varphi_A(Y))_{Y \in \mathcal{Y}_E}$ is a lattice isomorphism from $\mathcal{L}_{gr}(L_R(E))$ onto the product lattice $\mathcal{L}(R)^{\mathcal{Y}_E} = \prod_{Y \in \mathcal{Y}_E} \mathcal{L}(R)$, ordered componentwise.
     \end{enumerate}
     If $E$ is acyclic, then every ideal of $L_R(E)$ is graded by Theorem \ref{th:K}, so (2) and (3) apply to every ideal.
 \end{theorem}

\begin{proof}
By Lemma \ref{lem:structure}, $L_R(E) = \bigoplus_Y L_Y$. We first show that $A = \bigoplus_Y (A \cap L_Y)$ and that $A \cap L_Y$ is determined by $\varphi(Y)$; statements (1), (2) and (3) follow from this description.

Since the elements $\varepsilon_Y$ are central idempotents with sum $1$, every $a \in A$ can be written as $a = \sum_Y \varepsilon_Y a$, and each $\varepsilon_Y a$ belongs to $A \cap L_Y$, because $A$ and $L_Y$ are ideals and $\varepsilon_Y \in L_Y$. Hence $A = \bigoplus_Y (A \cap L_Y)$, where each $A \cap L_Y$ is an ideal of the ring $L_Y$. We now identify these ideals. Recall that $\varphi(u) = \{ k \ : \ ku \in A\}$ (Theorem \ref{th:iso}(2)).

Let $w$ be a sink. By Lemma \ref{lem:structure}(3), $L_{\{w\}}$ is isomorphic to $M_{n(w)}(R)$, and the vertex $w = ww^*$ corresponds to the diagonal matrix unit $E_{w,w}$. By (F4), the ideal $A \cap L_{\{w\}}$ corresponds to $M_{n(w)}(J)$, where $J$ is the set of scalars $k$ such that $kE_{w,w}$ corresponds to an element of $A$, that is, such that $kw \in A$. So $J = \varphi(w)$.

Now let $C$ be a cycle. By Lemma \ref{lem:structure}(4), $L_{C^0}$ is isomorphic to $M_{m(C)}(R[x,x^{-1}])$, and, by (F4) applied to the ring $R[x,x^{-1}]$, $A \cap L_{C^0}$ corresponds to $M_{m(C)}(J')$, where $J'$ is the ideal of $R[x,x^{-1}]$ formed by the Laurent polynomials $\sum_i k_ix^i$ such that $\sum_i k_ix^iE_{v_C,v_C}$ corresponds to an element of $A$, that is, such that $\sum_i k_iC^{\langle i \rangle} \in A$. Since $A$ is graded and the elements $C^{\langle i \rangle}$ have pairwise distinct degrees $i|C|$, the element $\sum_i k_iC^{\langle i \rangle}$ lies in $A$ if and only if each $k_iC^{\langle i \rangle}$ does. Multiplying by $C^{\langle -i \rangle}$, and using that $C^{\langle -i \rangle}C^{\langle i \rangle} = v_C$, we see that this happens if and only if $k_iv_C \in A$, that is, $k_i \in \varphi(C^0)$. Therefore $J' = \varphi(C^0)[x,x^{-1}]$.

This description gives (3) and the second isomorphism in (2). The map in (3) clearly preserves the order, and it is injective, because we have just seen that $A$ is determined by the family $(\varphi_A(Y))_Y$. For surjectivity, let $(J_Y)_Y$ be a family of ideals of $R$, and let $A'$ be the ideal of $L_R(E)$ that corresponds, through Lemma \ref{lem:structure}, to $\bigoplus_w M_{n(w)}(J_{\{w\}}) \oplus \bigoplus_C M_{m(C)}(J_{C^0}[x,x^{-1}])$. It is generated by the homogeneous elements $jw$, $j \in J_{\{w\}}$, and $jv_C$, $j \in J_{C^0}$, so it is graded, and the computation above gives $\varphi_{A'}(Y) = J_Y$. Thus the map in (3) is a bijection that preserves the order in both directions, that is, a lattice isomorphism. As for the quotient, $L_R(E)/A \cong \bigoplus_Y L_Y/(A \cap L_Y)$, and, by (F4), the summands are $M_{n(w)}(R/\varphi(w))$ and $M_{m(C)}\big(R[x,x^{-1}]/\varphi(C^0)[x,x^{-1}]\big) \cong M_{m(C)}\big((R/\varphi(C^0))[x,x^{-1}]\big)$.

For (1), fix $u \in E^0$ and $k \in R$. Since $u = \sum_Y \varepsilon_Yu$ and $A = \bigoplus_Y (A \cap L_Y)$, we have $ku \in A$ if and only if $k\varepsilon_Yu \in A$ for every cluster $Y$. By (\ref{eq:epsu}), $\varepsilon_Yu = \sum_{p \in P_Y, \, s(p) = u} pp^*$, which, in the matrix description of $L_Y$, is the diagonal matrix with entry $1$ at the positions $p \in P_Y$ that start at $u$ and $0$ elsewhere. If $u$ does not reach $Y$, there are no such paths and $\varepsilon_Yu = 0$. Otherwise this matrix is non-zero, and $k\varepsilon_Yu$ lies in $A \cap L_Y$, that is, it has its entries in $\varphi(Y)$, respectively in $\varphi(Y)[x,x^{-1}]$, if and only if $k \in \varphi(Y)$. Hence $ku \in A$ if and only if $k \in \varphi(Y)$ for every cluster $Y$ that $u$ reaches, which is (1).

It remains to prove the first isomorphism in (2). Since $\varphi$ is constant on clusters, each $X_I$ is a union of clusters. Comparing with the second isomorphism, it suffices to show that, for $I \neq R$ and $X = X_I$, the algebra $L_{R/I}(_XE)$ is the direct sum of the blocks $M_{n(w)}(R/I)$, for the sinks $w \in X$, and $M_{m(C)}\big((R/I)[x,x^{-1}]\big)$, for the cycles $C$ with $C^0 \subseteq X$. To this end, we show that $_XE$ is the disjoint union of the graphs $_YE$, $Y \subseteq X$, and that each $_YE$ has the same matrix block as $Y$ has in $E$. Recall from Definition \ref{Gporcupine} that $F(X)$ is the set of paths that enter $X$ only at their last vertex, and let $\alpha \in F(X)$. If a vertex of $\alpha$ other than its range were terminal, then, by (T1), the rest of $\alpha$ would stay in the cluster of that vertex; this cluster would contain $r(\alpha) \in X$, so it would be contained in $X$, against the definition of $F(X)$. Thus all the vertices of $\alpha$ except its range are non-terminal, and $\alpha \in F(Y)$ for the cluster $Y$ of $r(\alpha)$. Moreover, the edges of $E$ with both ends in $X$ are the edges of the cycles in $X$, because sinks emit no edges and the vertices of a cycle only emit the edge of the cycle. Consequently, $_XE$ is the disjoint union of the graphs $_YE$, $Y \subseteq X$. In $_YE$, every new vertex $w^\alpha$ emits exactly one edge, which takes it one step closer to $Y$. So $_{\{w\}}E$ is a finite acyclic graph whose only sink is $w$, and $_{C^0}E$ is a finite graph without sinks whose only cycle is $C$, which has no exits. Both graphs satisfy the hypotheses of Lemma \ref{lem:structure} and have a single cluster, so $L_{R/I}(_YE)$ is a matrix algebra over $R/I$, respectively $(R/I)[x,x^{-1}]$, whose size is the number of paths of $_YE$ ending at $w$, respectively ending at $v_C$ and not of the form $qC$. Each vertex of $_YE$ is the source of exactly one such path, obtained by following the unique edges towards $Y$, and then along $C$ until $v_C$ in the case of a cycle. So the size is the number of vertices of $_YE$, which is $|Y| + |F(Y)|$. This is also the size of the block of $Y$ in $E$: clearly $P_w = \{w\} \sqcup F(\{w\})$, and a path of $P_C$ either runs along $C$, and is then determined by its source in $C^0$, or starts with some $\alpha \in F(C^0)$ and then follows $C$ from $r(\alpha)$ to $v_C$ without completing a turn. Hence $L_{R/I}(_{\{w\}}E) \cong M_{n(w)}(R/I)$ and $L_{R/I}(_{C^0}E) \cong M_{m(C)}\big((R/I)[x,x^{-1}]\big)$, and summing over the clusters $Y \subseteq X_I$ and over $I \neq R$ gives the first isomorphism in (2).
\end{proof}

\begin{example}\label{ex:recovering} Let $E$ be the graph of Example \ref{ex:porcurpine1}, $R = \Z$, and $A$ the graded ideal with $\varphi(v_2) = p\Z$ and $\varphi(v_3) = q\Z$, for primes $p\neq q$; it is generated by $pv_2$ and $qv_3$. By Theorem \ref{Th:MainI}, $\varphi(u)= \varphi(v_1) = pq\Z$, and
$$L_\Z(E)/A \cong L_{\Z_p}\left(_{\{v_2\}}E\right) \oplus L_{\Z_q}\left(_{\{v_3\}}E\right) \cong M_3(\Z_p) \oplus M_3(\Z_q).$$
\end{example}

\begin{example}\label{ex:cycle}
Let $E$ be the graph
\[
\xymatrix{
 & & & \bullet_{v_3} \ar[dd]^{c_3} \\
\bullet_{v} & \bullet_{u} \ar[l]_{e_1} \ar[r]^{e_2} & \bullet_{v_2} \ar[ur]^{c_2} & \\
 & & & \bullet_{v_1} \ar[ul]^{c_1}
}
\]
with the cycle $C = c_1c_2c_3$, based at $v_C = v_1$. The clusters are $\{v\}$ and $C^0 =\{v_1,v_2,v_3\}$, $P_v = \{v, e_1\}$ and $P_C = \{v_1, c_3, c_2c_3, e_2c_2c_3\}$, so $n(v) = 2$ and $m(C) = 4 = |C^0| + |F(C^0)|$. Let $A$ be the graded ideal with $\varphi(v) = J_1$ and $\varphi(C^0) = J_2$; then, by Theorem \ref{Th:MainI}(1), $\varphi(u) = J_1 \cap J_2$, and, by Theorem \ref{Th:MainI}(2),
$$L_R(E)/A \cong L_{R/J_1}\left(_{\{v\}}E\right) \oplus L_{R/J_2}\left(_{C^0}E\right) \cong M_2(R/J_1) \oplus M_4\big((R/J_2)[x,x^{-1}]\big),$$
where 
$$
    _{\{v\}}E: \quad  \xymatrix{ && \\
\bullet_{v} & \bullet_{w^{e_1}} \ar[l]_{f^{e_1}} \\
& &
} \qquad _{C^0}E: \quad 
\xymatrix{
 & & & \bullet_{v_3} \ar[dd]^{c_3} \\ & \bullet_{w^{e_2}}  \ar[r]^{f^{e_2}} & \bullet_{v_2} \ar[ur]^{c_2} & \\
 & & & \bullet_{v_1} \ar[ul]^{c_1}
}
$$
For instance, if $R = \Z$, $J_1 = 3\Z$ and $J_2 = 5\Z$, this is $M_2(\Z_3) \oplus M_4(\Z_5[x,x^{-1}])$. If $J_1 = J_2 = I \neq R$, then $A$ is $I$-basic, and the result agrees with Theorem \ref{cor:main}.
\end{example}

Both hypotheses of Theorem \ref{Th:MainI} are needed. The finiteness of $E$ cannot be removed, even for row-finite acyclic graphs.

\begin{example}\label{ex:tree}
Let $E$ be the infinite binary tree
$$
    E: \quad
    \xymatrix@C=1.1pc@R=1.3pc{
     & & & \bullet_{v_0} \ar[dll] \ar[drr] & & & \\
     & \bullet \ar[dl] \ar[dr] & & & & \bullet \ar[dl] \ar[dr] & \\
     \bullet & & \bullet & & \bullet & & \bullet \\
     \vdots & & \vdots & & \vdots & & \vdots
    }
$$
with root $v_0$, in which every vertex emits exactly two edges, each of them to a new vertex. The graph $E$ is row-finite and acyclic, and it has no terminal vertices in the sense of \cite{VasPQ}: there are neither sinks nor cycles, and no infinite path $\beta$ is terminal \cite[Definition 5.3]{VasPQ}, because $T(\beta^0)$ contains the vertices of the branches leaving $\beta$, which do not connect to $\beta$. So $E$ has no clusters, and the analogue of Theorem \ref{Th:MainI}(3) would give a lattice of graded ideals with a single element. However, $L_R(E)$ has many graded ideals. For every vertex $w$, the set $T(w)$ is hereditary and saturated, since a vertex whose two children lie in $T(w)$ lies in $T(w)$. Hence the basic graded ideals $I(T(w))$, $w \in E^0$, are pairwise distinct, and $L_R(E)$ has infinitely many graded ideals, even when $R$ is a field.
\end{example}

The hypothesis that no cycle has an exit cannot be removed either.

\begin{example}\label{ex:Toeplitz2}
Consider again the Toeplitz graph $T$ of Example \ref{ex:ToeplitzGF}
\[ T:
\xymatrix{
\bullet_{u}\ar@(ul,dl)_{c} \ar@/_.1pc/[r]_e &
\bullet_{v}}
\]
and let $A$ be the graded ideal of $L = L_\Z(T)$ with $\varphi(v) = 2\Z$ and $\varphi(u) = 4\Z$, that is, the ideal generated by $2v$ and $4u$. If Theorem \ref{Th:MainI} extended to $T$, we would expect $L/A$ to be the direct sum of a block for the sink $v$ over $\Z/\varphi(v) = \Z_2$, with infinitely many rows indexed by the paths $v$ and $c^ke$, $k \geq 0$, that end at $v$, and a block for the loop $c$ over $(\Z/\varphi(u))[x,x^{-1}]$; that is, $L/A \cong M_\infty(\Z_2) \oplus \Z_4[x,x^{-1}]$. However, $B = L/A$ is indecomposable.

Suppose that $B = B_1 \oplus B_2$ with $B_1$ and $B_2$ non-zero ideals. We reduce modulo $2$ and use the Cuntz--Krieger Uniqueness Theorem over $\Z_2$. Since the generators $2v$ and $4u$ of $A$ lie in $2L$, we have $A \subseteq 2L$, and so $B/2B \cong L/2L$, which is $L_{\Z_2}(T)$ by (F3). This is the Leavitt path algebra of $T$ over the field $\Z_2$, and $T$ satisfies Condition (L), because its only cycle, the loop $c$, has the exit $e$. If $Q$ is a non-zero ideal of $L_{\Z_2}(T)$, the quotient map $L_{\Z_2}(T) \rightarrow L_{\Z_2}(T)/Q$ is not injective, so, by the Cuntz--Krieger Uniqueness Theorem, it vanishes on some vertex; thus $Q$ contains $u$ or $v$, and in both cases it contains $v = e^*ue$. Hence any two non-zero ideals of $L_{\Z_2}(T)$ intersect, and $L_{\Z_2}(T)$ is indecomposable. On the other hand, $B/2B \cong B_1/2B_1 \oplus B_2/2B_2$, and both summands are non-zero. Indeed, $4L \subseteq A$, because $4u$ and $4v = e^*(4u)e$ lie in $A$ and $L$ is spanned by the elements $\alpha\beta^* = \alpha\,r(\alpha)\,\beta^*$; so $4B = 0$, and if $B_i = 2B_i$, then $B_i = 4B_i = 0$. This contradicts the indecomposability of $L_{\Z_2}(T)$.
\end{example}

\subsection*{Acknowledgements}
The authors would like to thank K.M. Rangaswamy for his helpful feedback on the paper. The second and third authors would like to thank the host institution Istanbul Technical University, Department of Mathematics Engineering for the warm and productive environment during their research visits.

\section*{Declarations}

\subsection*{Funding}
The third author was supported by Proyecto puente
URJC M3031 (URJC Grant) “Lie algebras, Jordan systems and other nonassociative structures”, Proyecto
puente URJC M3329 (URJC Grant) “Lie algebras, Jordan systems and related structures”, Programa
propio de movilidad URJC 2023, and partially supported by Junta de Andalucía FQM264.

\subsection*{Conflicts of interests/Competing interests}
The authors declare that they have no known competing financial interests or personal relationships that could have appeared to influence the work reported in this paper.

\end{document}